\documentclass[a4paper,10pt]{llncs}
\usepackage[latin1]{inputenc}
\usepackage{amsmath,amssymb,amsfonts,bm}
\usepackage{dsfont}                          
 \let\mathbb\mathds     
\DeclareMathAlphabet\PazoBB{U}{fplmbb}{m}{n} 
\usepackage[english,frenchb]{babel}
\usepackage{cite}
\usepackage{graphicx,psfrag,subfigure}
\graphicspath{{./figures/}}

\newcommand \Fcal {\mathcal{F}}

\newcommand \Gfrak {\mathfrak{G}}
\newcommand \Tfrak {\mathfrak{T}}
\newcommand \Xfrak {\mathfrak{X}}

\newcommand   \Xset  {\mathbb{X}}
\newcommand   \XX    {\Xset}
\newcommand   \RR    {\mathbb{R}}

\newcommand \EE   {\mathsf{E}}

\DeclareMathOperator* \argmax {argmax}

\newcommand \ddiff {\mathrm{d}}

\newcommand \normPDF \varphi
\newcommand \normCDF \Phi

\author{Romain Benassi \and Julien Bect \and Emmanuel Vazquez}

\institute{SUPELEC, Gif-sur-Yvette, France}

\title{Bayesian optimization\\ using sequential Monte Carlo}

\begin{document}
\maketitle
\selectlanguage{english}
\begin{abstract}
  We consider the problem of optimizing a real-valued continuous
  function $f$ using a Bayesian approach, where the evaluations of $f$
  are chosen sequentially by combining prior information about $f$,
  which is described by a random process model, and past evaluation
  results. The main difficulty with this approach is to be able to
  compute the posterior distributions of quantities of interest
  which are used to choose evaluation points. In this article, we decide
  to use a Sequential Monte Carlo (SMC) approach.
\end{abstract}

\section{Overview of the contribution proposed}
\vspace{-1em}

We consider the problem of finding the global maxima of a function
$f:\XX\to\RR$, where $\XX\subset\RR^d$ is assumed bounded, using the
\emph{expected improvement}~(EI) criterion \cite{mockus:78:abmse,
  jones:98:ego}. Many examples in the literature show that the EI
algorithm is particularly interesting for dealing with the
optimization of functions which are expensive to evaluate, as is often
the case in design and analysis of computer experiments
\cite{santner}. However, going from the general framework expressed in
\cite{mockus:78:abmse} to an actual computer implementation is a
difficult issue.

The main idea of an EI-based algorithm is a Bayesian one: $f$ is
viewed as a sample path of a random process $\xi$ defined on
$\RR^d$. For the sake of tractability, it is generally assumed that
$\xi$ has a Gaussian process distribution conditionally to a parameter
$\theta \in \Theta \subseteq \RR^s$, which tunes the mean and
covariance functions of the process. Then, given a prior
distribution~$\pi_0$ on~$\theta$ and some initial evaluation results
$\xi(X_1), \ldots, \xi(X_{n_0})$ at $X_1, \ldots, X_{n_0}$, an (idealized)
EI algorithm constructs a sequence of evaluations points
$X_{n_0+1}, X_{n_0+2},\ldots$ such that, for each $n \geq n_0$,
\begin{equation}
  \label{eq:fully_Bayesian_EI}
  X_{n+1} = \argmax_{x\in\XX} \bar\rho_n := \int_{\theta\in\Theta} \rho_n(x;\theta) \ddiff\pi_n(\theta)\,,
\end{equation}
where $\pi_n$ stands for the posterior distribution of $\theta$,
conditional on the $\sigma$-algebra $\Fcal_n$ generated by $X_1,
\xi(X_1), \ldots, X_n, \xi(X_n)$, and
\begin{equation*}
  \rho_n(x;\theta) := %
  \EE_{n,\theta}((\xi(X_{n+1}) - M_n)_+ \mid X_{n+1} = x) 
\end{equation*}
is the EI at $x$ given $\theta$, with $M_n = \xi(X_0)\vee\cdots\vee
\xi(X_n)$ and $\EE_{n,\theta}$ the conditional expectation given
$\Fcal_n$ and $\theta$. In practice, the computation of $\rho_n$ is
easily carried out (see \cite{jones:98:ego}) but the answers to the
following two questions will probably have a direct impact on the
performance and applicability of a particular implementation:
a)~How~to deal with the integral in $\bar\rho_n$?  b)~How~to deal with
the maximization of $\bar\rho_n$ at each step?

We can safely say that most implementations---including the popular
EGO algorithm \cite{jones:98:ego}---deal with the first issue by using
an \emph{empirical Bayes} (or \emph{plug-in}) approach, which consists
in approximating~$\pi_n$ by a Dirac mass at the maximum likelihood
estimate of $\theta$. A plug-in approach using maximum a posteriori
estimation has been used in~\cite{lizotte}; \emph{fully Bayesian}
methods are more difficult to implement (see \cite{lion5} and
references therein).  Regarding the optimization of $\bar \rho_n$ at
each step, several strategies have been proposed (see, e.g.,
\cite{jones:98:ego, diceOptim, bardenet:2010, gramacy:2011}).


This article addresses both questions simultaneously, using a
sequential Monte Carlo (SMC) approach \cite{chopin:2002,
  delmoral:2006:smc} and taking particular care to control the
numerical complexity of the algorithm. The main ideas are the
following.  First, as in~\cite{gramacy:2011}, a weighted sample
$\Tfrak_n=\{(\theta_{n,i}, w_{n,i})\in\Theta\times\RR, 1\leq i\leq
I\}$ from~$\pi_n$ is used to approximate~$\bar\rho_n$;~that is,
$\sum_{i=1}^I w_{n,i} \;\rho_n(x;\theta_{n,i}) \to_I \bar\rho_n(x)
$. Besides, at each step~$n$, we attach to each $\theta_{n,i}$ a
(small) population of candidate evaluation points $\{x_{n,i,j}, 1\leq
j \leq J \}$ which is expected to cover promising regions for that
particular value of $\theta$ and such that $\max_{i,j}
\bar\rho_n\left( x_{n,i,j} \right) \approx \max_x \bar\rho_n(x)$. 



\vspace{-1em}

\section{Algorithm and results}
\vspace{-1em}

At each step $n \geq n_0$ of the algorithm, our objective is to
construct a set of weighted particles
\begin{align}
  \label{equ:Gn}
  \Gfrak_{n} = \big\{\; & \left( \gamma_{n,i,j},
    w^\prime_{n,i,j} \right)\,,~ \nonumber \\ 
  & \;\, \gamma_{n,i,j} = (\theta_{n,i}, x_{n,i,j}) \in \Theta\times \XX,
  w^\prime_{n,i,j}\in\RR\,,~ 1 \le i \le I,\, 1 \le j \le J \;\big\}  
\end{align}
so that $\sum_{i,j} w^\prime_{n,i,j} \delta_{\gamma_{n,i,j}} \to_{I,J}
\pi'_n$, with 
$$
\ddiff\pi^{\prime}_n(\gamma) \;=\; \tilde g_n(x\mid\theta)\, \ddiff\lambda(x)\,\ddiff
\pi_n(\theta) \,,\quad x\in\XX\,,~\theta\in\Theta\,,~\gamma=(\theta,x),
$$ 
where $\lambda$ denotes the Lebesgue measure, $\tilde g_n(x\mid\theta)
= g_n(x\mid\theta)/c_n(\theta)$, $g_n(x\mid\theta)$ is a criterion
that reflects the interest of evaluating at $x$ (given $\theta$ and
past evaluation results), and $c_n(\theta) = \int_{\XX}
g_n(x\mid\theta) \ddiff x$ is a normalizing term. For instance, a
relevant choice for $g_n$ is to consider the probability that $\xi$
exceeds $M_{n}$ at $x$, at step $n$. (Note that we consider less
$\theta$s than $x$s in $\Gfrak_{n}$ to keep the numerical complexity
of the algorithm low.)

To initialize the algorithm, generate a weighted sample
$\Tfrak_{n_0}=\{ (\theta_{n_0,i}, w_{n_0,i}),$ $1\leq i\leq I\}$ from
the distribution $\pi_{n_0}$, using for instance importance sampling
with $\pi_0$ as the instrumental distribution, and pick a density
$q_{n_0}$ over $\XX$ (the uniform density, for example). 
Then, for each $n \geq n_0$: %
\\[0.6em]
\emph{Step 1: demarginalize} --- Using~$\Tfrak_n$ and~$q_n$, construct
a weighted sample~$\Gfrak_n$ of the form~\eqref{equ:Gn}, with
$x_{n,i,j} \stackrel{\text{\tiny iid}}{\sim} q_n$, $ w^\prime_{n,i,j}
= w_{n,i} \frac{g_n(x_{n,i,j}| \theta_{n,i}) }{ q_n(x_{n,i,j})
  c_{n,i}}$, and $c_{n,i} =\sum_{j'=1}^{J} \frac{g_n(x_{n,i,j'}|
  \theta_{n,i})}{ q_n(x_{n,i,j'}) }$. %
\\[0.6em]
\emph{Step 2: evaluate} --- Evaluate~$\xi$ at $X_{n+1} = \argmax_{i,j}
\sum_{i'=1}^I w_{n,i'}\, \rho_n( x_{n,i,j}; \theta_{n,i'} )$. %
\\[0.6em]
\emph{Step 3: reweight/resample/move} --- Construct~$\Tfrak_{n+1}$
from~$\Tfrak_n$ as in \cite{chopin:2002}: reweight the
$\theta_{n,i}$s using $ w_{n+1,i} \propto \frac{
  \pi_{n+1}(\theta_{n,i}) }{ \pi_{n}(\theta_{n,i}) } \, w_{n,i}$,
resample (e.g., by multinomial resampling), and move the
$\theta_{n,i}$s to get $\theta_{n+1,i}$s using an independant
Metropolis-Hastings kernel.
\\[0.6em]
\emph{Step 4: forge $q_{n+1}$} --- Form an estimate~$q_{n+1}$ of the
second marginal of~$\pi'_n$ from the weighted sample $\Xfrak_n = \{
(x_{n,i,j}, w'_{n,i,j}), 1\le i \le I, 1 \le j \le J \}$. Hopefully,
such a choice of~$q_{n+1}$ will provide a good instrumental density
for the next demarginalization step. Any (parametric or
non-parametric) density estimator can be used, as long as it is easy
to sample from; in this paper, a tree-based histogram estimator is
used.

\emph{Nota bene:} when possible, some components of~$\theta$ are
integrated out analytically in~\eqref{eq:fully_Bayesian_EI} instead of
being sampled from; see~\cite{lion5}.

\textbf{Experiments.} Preliminary numerical results, showing the
relevance of a fully Bayesian approach with respect to empirical Bayes
approach, have been provided in \cite{lion5}. The scope of these
results, however, was limited by a rather simplistic implementation
(involving a quadrature approximation for~$\bar\rho_n$ and a
non-adaptive grid-based optimization for the choice of~$X_{n+1}$). We
present here some results that demonstrate the capability of our new
SMC-based algorithm to overcome these limitations.

The experimental setup is as follows. We compare our SMC-based
algorithm, with $I=J=100$, to an EI algorithm in which: 1) we fix
$\theta$ (at a ``good'' value obtained using maximum likelihood
estimation on a large dataset); 2) $X_{n+1}$ is obtained by exhaustive
search on a fixed LHS of size $I\times J$. In both cases, we consider
a Gaussian process $\xi$ with a constant but unknown mean function
(with a uniform distribution on $\RR$) and an anisotropic Mat\'ern
covariance function with regularity parameter $\nu=5/2$.  Moreover,
for the SMC approach, the variance parameter of the Mat\'{e}rn
covariance function is integrated out using a Jeffreys prior and the
range parameters are endowed with independent lognormal priors.

\begin{figure}[h!]
\begin{minipage}{0.48\linewidth}
  \psfrag{ylabel}[t][t]{\scriptsize $\log(\max f -M_n)$}
\psfrag{iterations}[t][t]{\scriptsize number of function evaluations}
\psfrag{distance au maximum : 100}{}
\psfrag{ego}{\scriptsize ref EI}
\psfrag{smc Bayes}{\scriptsize EI+SMC}
\centering\subfigure[Branin function (dimension 2)]
{\label{fig:branin}\includegraphics[width=6cm]{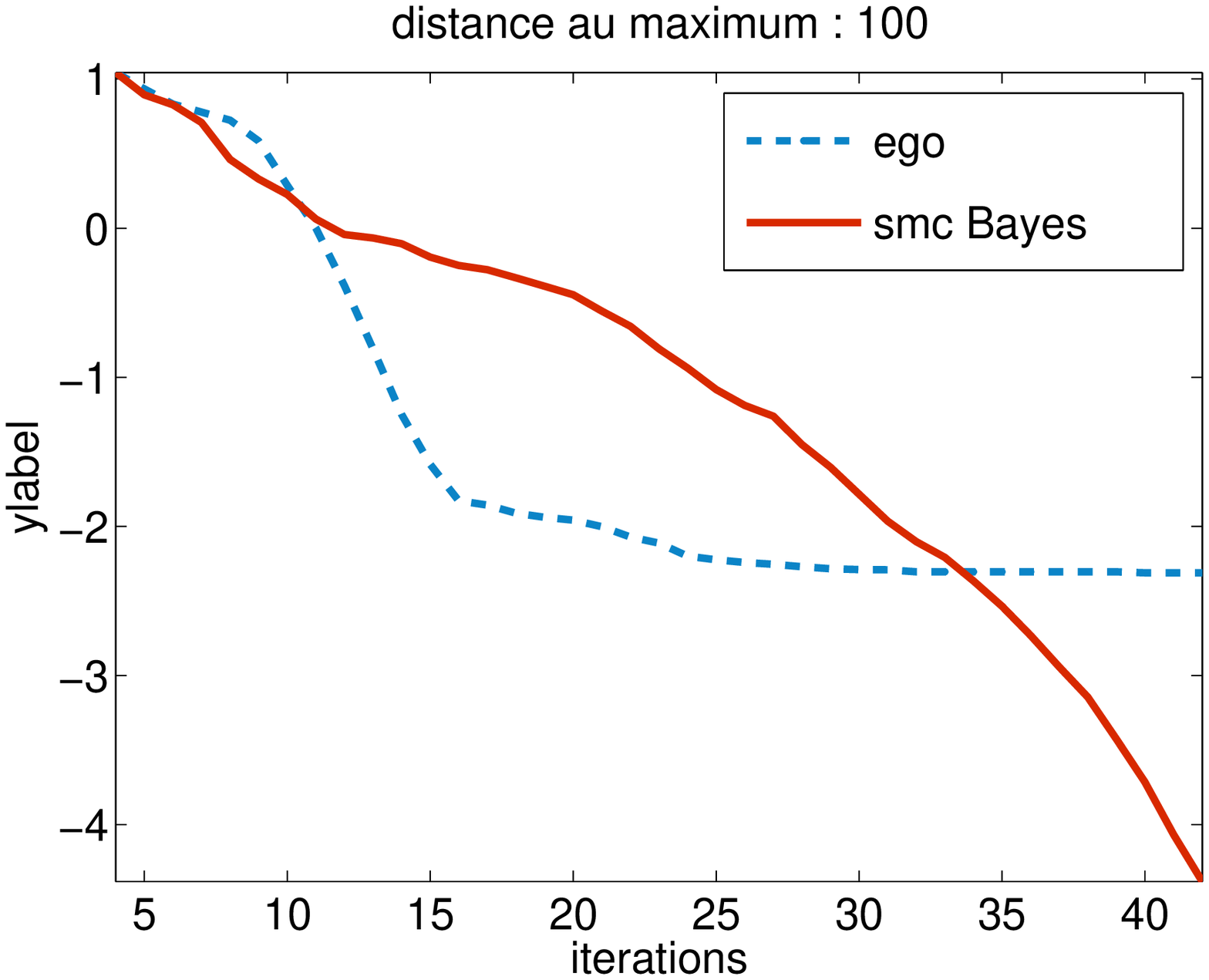}}
\end{minipage}
\hfill
\begin{minipage}{0.48\linewidth}
  \psfrag{ylabel}[t][t]{\scriptsize $\log(\max f -M_n)$}
\psfrag{iterations}[t][t]{\scriptsize number of function evaluations}
\psfrag{distance au maximum : 100}{}
\psfrag{ego}{\scriptsize ref EI}
\psfrag{smc Bayes}{\scriptsize EI+SMC}
\centering\subfigure[Hartmann 6 function (dimension 6)]
{\label{fig:hart6}\includegraphics[width=6cm]{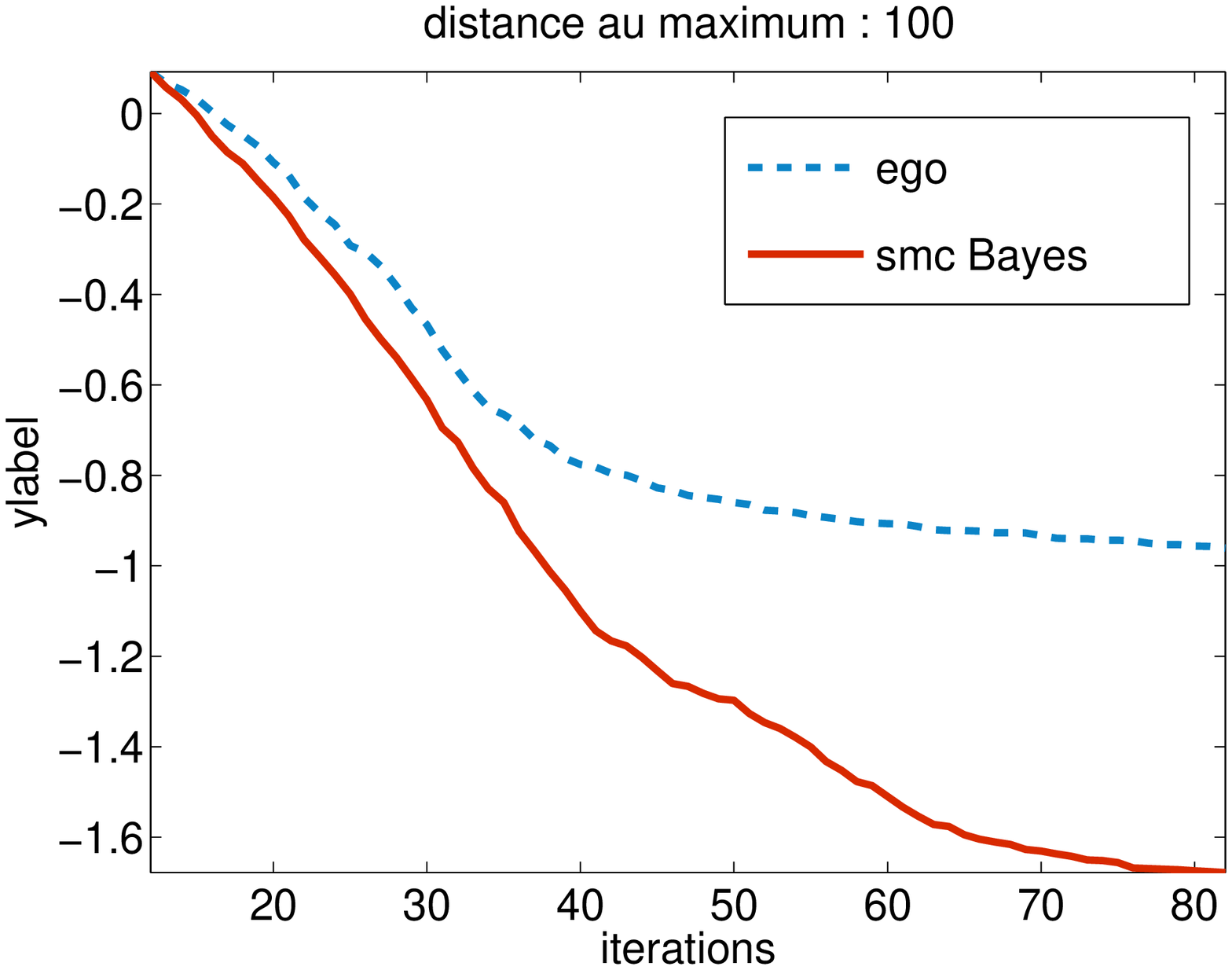}}
\end{minipage}
  \caption{A comparison of the average error to the maximum (100 runs)}
\label{fig:average_error}
\end{figure}

\textbf{Results.} Figures~\ref{fig:branin} and~\ref{fig:hart6} show
the average error over $100$ runs of both algorithms, for the Branin
function ($d=2$) and the log-transformed Hartmann~6 function ($d =
6$). For the Branin function, the reference algorithm performs better
on the first iterations, probably thanks to the ``hand-tuned''
parameters, but soon stalls due to its non-adaptive search strategy.
Our SMC-based algorithm, however, quickly catches up and eventually
overtakes the reference algorithm. On the Hartmann~6 function, we
observe that the reference algorithm always lags behind our new
algorithm. 

We have been able to find results of this kind for other test
functions. These findings are promising and need to be further
investigated in a more systematic large-scale benchmark study.



\bibliographystyle{plain}

\end{document}